\documentclass[11pt]{article}
\usepackage{graphicx}
\usepackage{float}
\usepackage{amsmath}
\usepackage{amsfonts}
\usepackage{paralist}
\usepackage{amsfonts}
\usepackage{enumerate}
\usepackage{amssymb}
\usepackage{graphicx}
\usepackage{amsthm}
\usepackage{a4wide}
\usepackage{array}
\usepackage{listings}
\usepackage{color}
\usepackage{tikz}
\usepackage{pgfplots}
\usepackage{mathtools}
\usetikzlibrary{arrows,positioning,calc}
\newcommand {\mat}[1]{\left[\begin{array}{#1}}
\newcommand {\rix}          {\end{array}\right]}
\newcommand {\eq}       [1] {\begin{equation}\label{#1}}
\newcommand {\en}           {\end{equation}}
\newcommand {\mpar} [1] {\marginpar{\fussy\tiny #1}}
\renewcommand {\mpar} [1] {}
\newtheorem{theorem}          {Theorem}[section]
\newtheorem{lemma}         [theorem]{Lemma}
\newtheorem{corollary}     [theorem]{Corollary}

\newcommand{\bal}{\begin{center}\begin{tabular}{ll}}
\newcommand{\eal}{\end{tabular}\end{center}}
\newcommand{\itxt}[1]{\noindent The following table lists the input
  arguments of the #1. \\}
\newcommand{\otxt}[1]{\noindent The following table lists the output
  arguments of the #1. \\}
\newcommand{\C} {{\mathbb C}}
\newcommand{\R} {{\mathbb R}}
\newcommand{\LL} {{\mathbb L}}
\newcommand {\rank}     {\mathop{\rm rank}\nolimits}
\catcode`@=11
\font\tenex=cmex10 
\newdimen\p@renwd
\setbox0=\hbox{\tenex B} \p@renwd=\wd0 
\def\bmat#1{\begingroup \m@th
  \setbox\z@\vbox{\def\cr{\crcr\noalign{\kern2\p@\global\let\cr\endline}}%
    \ialign{$##$\hfil\kern2\p@\kern\p@renwd&\thinspace\hfil$##$\hfil
      &&\quad\hfil$##$\hfil\crcr
      \omit\strut\hfil\crcr\noalign{\kern-\baselineskip}%
      #1\crcr\omit\strut\cr}}%
  \setbox\tw@\vbox{\unvcopy\z@\global\setbox\@ne\lastbox}%
  \setbox\tw@\hbox{\unhbox\@ne\unskip\global\setbox\@ne\lastbox}%
  \setbox\tw@\hbox{$\kern\wd\@ne\kern-\p@renwd\left[\kern-\wd\@ne
    \global\setbox\@ne\vbox{\box\@ne\kern2\p@}%
    \vcenter{\kern-\ht\@ne\unvbox\z@\kern-\baselineskip}\,\right]$}%
  \null\;\vbox{\kern\ht\@ne\box\tw@}\endgroup}
\catcode`@=12
\newtheorem{remark} [theorem] {Remark}

\begin{document}
\thispagestyle{empty}
\begin{center}

\renewcommand{\thefootnote}{\fnsymbol{footnote}}

\begin{minipage}{0.29\textwidth}
	\includegraphics[width=1in]{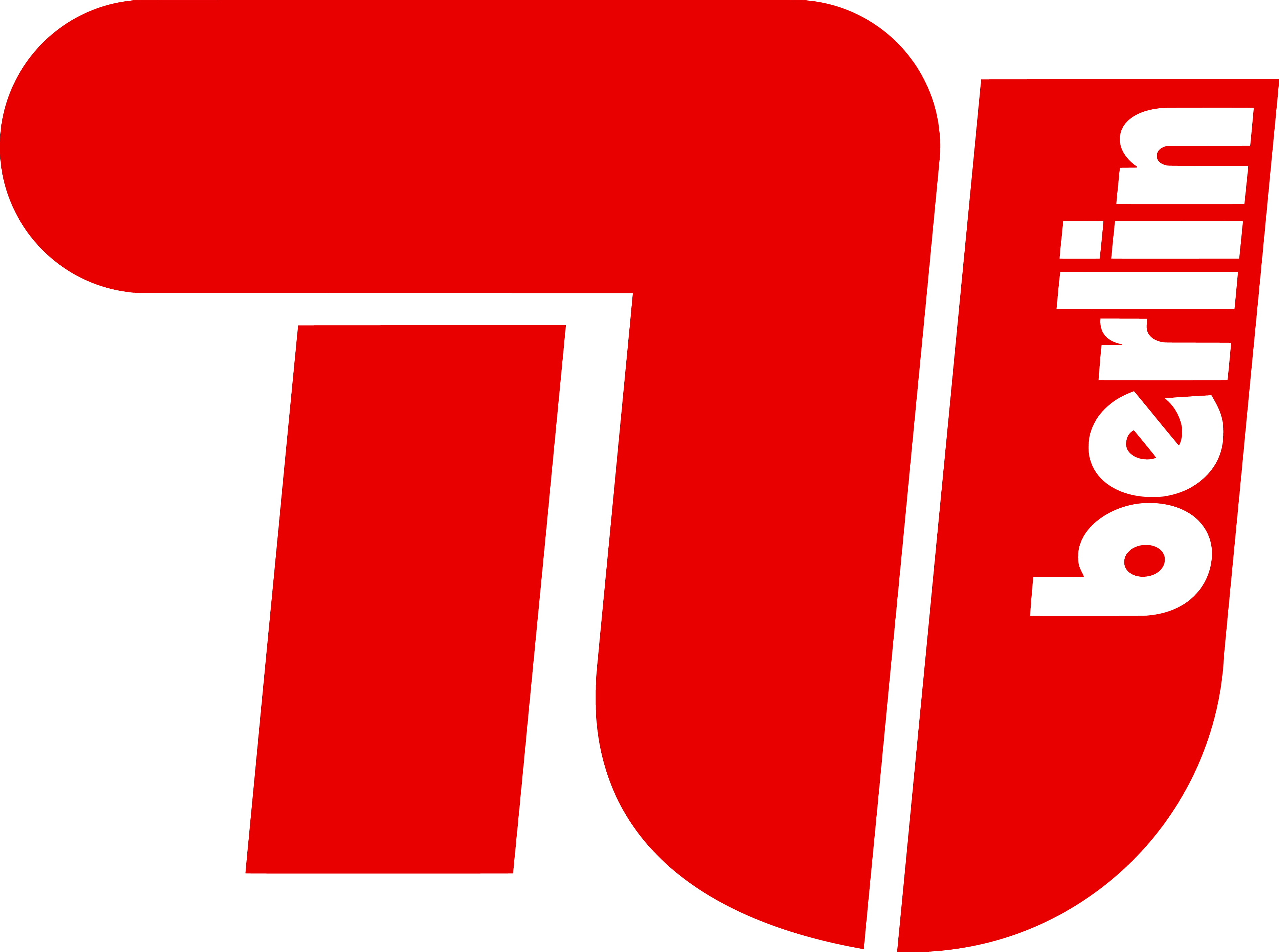}
\end{minipage}
\begin{minipage}{0.69\textwidth}
	\begin{flushright}
		{\Huge Technische Universit\"at Berlin\\}
		\Large Institut f\"ur Mathematik
	\end{flushright}
\end{minipage}
\\ [50mm]
\scalebox{0.95}{\bf\Large A \texttt{Matlab} Toolbox for the Regularization of Descriptor Systems}\\[0.2cm]
\scalebox{0.95}{\bf\Large Arising from Generalized Realization Procedures}\\[13mm]
{\large\bf A. Binder \qquad V. Mehrmann \qquad A. Miedlar \qquad P. Schulze} \\ [18mm]
{\bf Preprint 24-2015} \\
\vfill {\bf Preprint-Reihe des Instituts f\"ur Mathematik} \\[1mm]
{\bf Technische Universit\"at Berlin} \\[1mm]
{\bf \tt http://www.math.tu-berlin.de/preprints} \\[8mm]
{\bf Preprint 24-2015 \hfill December 2015}
\end{center}

\date{\today}
\title{A \texttt{Matlab} Toolbox for the Regularization of Descriptor Systems Arising from Generalized Realization Procedures}

\author{A. Binder\footnotemark[1] \and V. Mehrmann\footnotemark[1] \and A. Miedlar\footnotemark[2] \and P. Schulze\footnotemark[1]}

\maketitle
\renewcommand{\thefootnote}{\fnsymbol{footnote}}
\footnotetext[1]{
	Institut f\"ur Mathematik, TU Berlin, Germany, \texttt{$\{$binder,pschulze,mehrmann$\}$@math.tu-berlin.de}}
\footnotetext[2]{
University of Minnesota, Minneapolis, USA \texttt{amiedlar@umn.edu}.
}
\noindent In this report we introduce a \texttt{Matlab} toolbox for the regularization of descriptor systems. We apply it, in particular, for systems resulting from the generalized realization procedure of \cite{MayA07}, which generates, via rational interpolation techniques,
a linear \emph{descriptor system} from interpolation data. The resulting system needs to be regularized to make it feasible for the use in simulation, optimization, and control. This process is called \emph{regularization}.

\section{Descriptor Systems} \label{intro}
We follow the notation and the basic concepts of \cite{BunBMN99}.
A linear {\em descriptor system} is of the form
\begin{subequations}
\begin{align}
                E \dot{x} & = Ax + B u,\label{a-11}\\
                y         & = C x+  D u,\label{a-12}
                \end{align}
                \label{eq:a}
\end{subequations}
where
\(E,\,A \in \mathbb R^{n,n}\),
\(B \in \mathbb R^{n,m}\),
\(C \in \mathbb R^{p,n}\),
\(D \in \mathbb R^{p,m}\),
                and
\(\dot{x} = dx/dt\).
The response  of  a descriptor system can be  described in terms of
the eigenvalues of the matrix pencil \(\alpha E - \beta A\), which is
said to be {\em regular} if
\(\det(\alpha E - \beta A) \ne 0\) for some
\((\alpha,\,\beta) \in \C^2\).
For regular pencils, {\em generalized eigenvalues} are
the pairs \((\alpha,\beta) \in \C^2 \setminus \{(0,0)\}\),
for which \(\det(\alpha E - \beta A ) = 0\).
If \(\beta \neq 0\), then the pair represents the finite
eigenvalue \(\lambda = \alpha/\beta\).
If \(\beta = 0\), then \((\alpha,\beta)\) represents an
infinite eigenvalue.

In frequency domain, for zero initial conditions $x(t_0)=0$ and a regular pencil \(\alpha E - \beta A\), there exists the rational \emph{transfer function}
\begin{equation}\label{transfunction}
 H(s)=C(sE-A)^{-1}B+D,
\end{equation}
which maps Laplace transforms of the input functions $u$ to the Laplace transforms of the corresponding output functions $y$. A finite eigenvalue \(\lambda = \alpha/\beta\) is a pole
of the transfer function of the descriptor system
\eqref{eq:a}.

In the following we denote a matrix with
orthonormal columns spanning the right nullspace of the matrix $M$
by $S_{\infty}(M)$ and a matrix with orthonormal
columns spanning the left nullspace of $M$ by $T_{\infty}(M)$.
These matrices are not uniquely
determined although the  spaces are, but for ease of notation, we  speak of these matrices as the corresponding spaces.

For regular pencils the solution of the system equations can be
characterized in terms of the Weierstra\ss\   Canonical Form (WCF),
\cite{Gan59b}.
\begin{theorem} {\bf Weierstra\ss\  Canonical Form}
\label{r-wcf-theorem}
If \({\alpha E - \beta A}\) is a regular pencil,
then there exist nonsingular matrices
\(X = [ X_r,\,X_\infty ] \in \mathbb R^{n,n}\)
and
\(Y = [ Y_r,\,Y_\infty ] \in \mathbb R^{n,n}\)
for which
\begin{equation}\label{r-wcf-1}
                        Y^T E X =
                        \mat{c} Y_r^T \\* Y_\infty^T \rix
                        E
                        \mat{cc} X_r   &  X_\infty   \rix
                        =
                        \mat{cc} I & 0 \\* 0 & N \rix ,
\end{equation}
and
\begin{equation}\label{r-wcf-2}
                        Y^T A X =
                        \mat{c} Y_r^T \\* Y_\infty^T \rix
                        A
                        \mat{cc} X_r   &  X_\infty   \rix
                        =
                        \mat{cc} J & 0 \\* 0 & I \rix ,
\end{equation}
where $J$ is a matrix in Jordan canonical form
whose diagonal elements are the finite eigenvalues of
the pencil and $N$ is a nilpotent matrix, also in
Jordan form.
\(J\) and \(N\) are unique up to permutation of Jordan
blocks.
\end{theorem}
\noindent The {\em index} $\nu$ of the pencil $\alpha E-\beta A$ is the index of nilpotency of the
nilpotent matrix $N$ in (\ref{r-wcf-1}).
By convention, if \(E\) is nonsingular, the pencil is said to be
of index zero.
A descriptor system is regular and of index at most
one if and only if it has exactly $q = \rank(E)$ finite
eigenvalues. The following lemma of \cite{KauNC89} gives a useful characterization of regular,
index one pencils.
\begin{lemma}  \label{r-index-lemma}
                The following statements are equivalent:
\begin{enumerate}
\item
The pencil $\alpha E - \beta A$ is regular and has
index less than or equal to one.
\item
\(
\rank \left( \mat{c} E \\* T_\infty^T(E) A \rix \right ) =
\rank \left( E + T_\infty(E) T_\infty^T(E) A \right ) = n
\).
\item
\(
\rank \left( [ E,\, A S_\infty(E) ] \right) =
\rank \left( E + A S_\infty(E) S_\infty^T(E) \right )
= n
\).
\item
\(T_\infty (E)^T A S_\infty(E)\) is nonsingular.
\item
If
\[
U^T E V = \mat{cc} \Sigma_r & 0 \\ 0 & 0 \rix
\]
is the singular value decomposition (SVD) of \(E\) (with
orthogonal matrices \(U,\,V\) and a nonsingular,
diagonal matrix \(\Sigma_r\in\mathbb R^{r\times r}\)), then
the \((n-r)\times (n-r)\) bottom right matrix $A_{22}$ of \( U^TAV\)
is nonsingular.
\end{enumerate}
\end{lemma}
\noindent In the notation of (\ref{r-wcf-1})--(\ref{r-wcf-2}),
classical solutions of (\ref{a-11}) take the form
\[
                x(t) = X_r z_1(t) + X_\infty z_2(t),
\]
where
\begin{eqnarray*}
	\dot{z}_1   & = & J z_1 + Y_r^T      B u\\*
		\nonumber
	N \dot{z}_2 & = &   z_2 + Y_\infty^T B u
\end{eqnarray*}
and one has the explicit solution
        \begin{eqnarray}
                \nonumber z_1(t) & = &
                        e^{tJ} z_1(0)
                        + \int_0^t e^{(t-s)J} Y_r^TB u(s) \, ds , \\*
        \label{a-17}    z_2(t) & = & -\sum_{i=0}^{\nu-1} \frac{d^i}{dt^i}
                        \left(N^i Y_\infty^T B u(t) \right).
        \end{eqnarray}
Equation (\ref{a-17}) shows that the input functions must belong
to some suitable function space ${\cal U}_{ad}$ and, to ensure a smooth response for every continuous
input \(u(t)\), it is necessary for the system to be regular and
have index less than or equal to one. Moreover, the possible values of the
initial condition \(x(0)\) are restricted.
The initial state must be a member of the set of {\em consistent }
initial conditions, i.\,e.,
\[
\mathcal S
\equiv
\left\{ X_r z_1 + X_\infty z_2 \
\left | \
                                z_1 \in \R^r, \\  z_2 =
                                -\sum_{i=0}^{\nu-1}
                             \left (
                                \frac{d^i}{dt^i} ( N^i Y_\infty^T Bu)(0)
                             \right),
                        \right.  u(t)\in {\cal U}_{ad}
\right\}.
\]
The set of {\em reachable} states
of (\ref{a-11}) from the {\em solution space} set \(\mathcal S\) of consistent initial conditions is $\mathcal S$ itself.

\section{Realization\label{sec:realization}}
In \cite{MayA07} a method for the \emph{generalized realization problem} was presented. From given interpolation data, obtained by measurements from a real system or numerical simulation via a mathematical
model, it generates a \emph{descriptor  system} of the form \eqref{eq:a}, i.\,e.,
\begin{equation}
\begin{aligned}
E\dot x(t)&=&Ax(t)+Bu(t),\\
y(t)&=&C x(t)+Du(t).
\end{aligned}
\label{descsystem}
\end{equation}
%
%
The generalized realization problem of \cite{MayA07} deals mainly with two cases.
\begin{enumerate}
\item In the scalar interpolation case, the given data consist of a vector of interpolation points $s=[s_i]\in \mathbb C^{N}$ and a vector of interpolation values  $f=[f_i]\in \mathbb{C}^{N}$, and the realization problem constructs a transfer function of the form~(\ref{transfunction}) satisfying the interpolation conditions
\begin{equation*}
	H(s_i)=C(s_iE-A)^{-1}B+D=f_i \qquad i=1,\dots,N.
\end{equation*}
\item In the matrix interpolation case, the interpolation points are again contained in a vector $s=[s_i]\in \mathbb C^{N}$. However, the interpolation values are summarized in form of a block matrix $F=[F_i]$  of $N$ matrices $F_i\in \mathbb{C}^{p\times m}$ and the interpolation problem takes the following form:
First, \emph{right and left tangential data} are sampled by multiplying the matrix data $F_i$ from the right (left) with arbitrary right (left) tangential directions such that as \emph{right tangential interpolation conditions}  we get
\begin{equation}\label{rightintcond}
H(\lambda_i)r_i=(C(\lambda_i E-A)^{-1}B+D)r_i=w_i, \qquad i=1,\dots ,\rho,
\end{equation}
and as \emph{left tangential interpolation conditions}  we get
\begin{equation}\label{leftintcond}
	l_jH(\mu_j)=l_j(C(\mu_j E-A)^{-1}B+D)=v_j, \qquad j=1,\dots ,\nu,
\end{equation}
where $r_i$ ($l_j$) are the right (left) tangential directions, $w_i$ ($v_j$) are the right (left) tangential values, and $\lambda_i$ ($\mu_j$) are the right (left) interpolation points which are a subset of $\{s_1,\ldots,s_N\}$.
\end{enumerate}
The interpolation technique is realized in the \texttt{Matlab} codes \texttt{realization}
and \texttt{loewner\underline{~}mod} where the latter one is (a slightly modified version of) an m-File provided by the authors of \cite{MayA07}. Analytically, it can be shown that the obtained realization \eqref{descsystem} is regular and \emph{minimal} (and thus \emph{controllable and observable}), see \cite{MayA07}. However, there are no results regarding the index of the obtained descriptor system. Moreover, when the realization is computed numerically, the analytically guaranteed properties of regularity and minimality may be lost due to finite precision  arithmetic. Thus, in general the realization obtained by computation may be non-regular, have index larger than one and miss certain controllability and observability properties, and therefore requires a regularization procedure which is described in the next section.

\section{Controllability and Observability Conditions}
      \label{controllability}

Given the descriptor system (\ref{descsystem}), one or more of the following conditions are essential
for most classical design aims, see e.g. \cite{Ber13,BunBMN99,Dai89}.
\begin{equation}\label{a-20}
   \mbox{
   \begin{tabular}{ll}
     {\bf C0}:&
        \parbox[t]{3.5in} {
          \(\rank [ \alpha E - \beta A,\, B ] = n\)
        for all \( (\alpha,\beta)\in \C^2\backslash\{(0,0)\}\).
        } \\
     {\bf C1}:&
        \parbox[t]{3.5in} {
          \(\rank [ \lambda E - A,\, B ] = n\) for all \(\lambda \in \C\).
        } \\
     {\bf C2}:&
        \parbox[t]{3.5in} {
           \(\rank [ E,\, A S_\infty(E),\, B ] = n\).
        }
   \end{tabular}
   }
\end{equation}
A regular system is {\em completely controllable} or {\em C-controllable} if
{\bf C0} holds  and is
{\em strongly controllable} or  {\em S-controllable}
if {\bf C1} and {\bf C2} hold \cite{BunBMN99}.
Complete controllability ensures
that for any given initial and final states
\(x_0,\, x_f\)
there exists an admissible control that transfers
the system from $x_0$ to $x_f$ in finite time,
while strong controllability ensures the same for
any given initial and final states
\(x_0,\, x_f \in \mathcal S\) (the solution space).

Regular systems that satisfy
condition~{\bf C2} are called {\em  controllable at infinity\/}
or {\em impulse controllable\/} \cite{Dai89}.
For these systems, impulsive modes can be excluded by a suitable
linear feedback.\\

\noindent Observability for descriptor systems is the
dual of controllability. We define the following conditions:
\eq{a-22}
   \mbox{
   \begin{tabular}{ll}
     {\bf O0}:&
        \parbox[t]{3.5in}{
           \(\rank\mat{c}\alpha E-\beta A\\C\rix = n\) for all
\( (\alpha,\beta) \in \C^2\backslash\{(0,0)\}\).
        } \\*[0.3cm]
     {\bf O1}:&
        \parbox[t]{3.5in}{
           \(\rank\mat{c}\lambda E-A\\C\rix = n\) for all \(\lambda \in \C\).
        } \\*[0.3cm]
     {\bf O2}:&
        \parbox[t]{3.5in}{
           \(\rank \mat {c} E \\ T_{\infty}^T(E) A \\ C \rix = n\).
        }
   \end{tabular}
   }
\end{equation}
It is immediate that condition {\bf O0} implies {\bf O1} and {\bf O2}.
Moreover, {\bf O1}  and
\begin{equation}\label{a-23}
	\rank \mat {c}  E \\ C \rix  = n,
\end{equation}
together hold if and only if {\bf O0} holds.
A regular descriptor system is called
{\em completely observable} or {\em C-observable} if condition {\bf O0} holds and
is called  {\em strongly observable} or {\em S-observable}
if conditions {\bf O1} and {\bf O2} hold.
A regular system that satisfies condition {\bf O2} is called
{\em observable at infinity} or {\em impulse-observable}.

Conditions (\ref{a-20})--(\ref{a-23}) are preserved under non-singular
equivalence transformations as well as under state and output feedback,
i.\,e., if the system  satisfies {\bf C0},
{\bf C1}, or {\bf C2}, 
then for any non-singular
\(U \in \mathbb R^{n,n}\), \(V \in \mathbb R^{n,n}\),
\(W \in \mathbb R^{m,m}\)  and for any \(F_1 \in \mathbb R^{m,n}\) and \(F_2 \in \mathbb R^{m,p}\),
the system  $(\tilde E, \tilde A, \tilde B, \tilde C)$, where
\begin{equation}\label {a-25}
    \tilde E = U E V,  \qquad \tilde A = U A V, \qquad \tilde B = U B W
\end{equation}
or
\[
   \tilde E = E, \qquad \tilde A = A + B F_1 ,  \qquad \tilde B = B
\]
or
\[
   \tilde E = E, \qquad \tilde A = A + B F_2 C ,  \qquad \tilde B = B
\]
also satisfies these conditions.
Analogous properties hold for {\bf O0}, {\bf O1} and {\bf O2}.

\section{Regularization}
In general, due to the finite precision  arithmetic, it cannot be guaranteed that the system computed by the realization procedure presented in \cite{MayA07} satisfies the described regularity, controllability and observability conditions of Section~\ref{controllability}. Therefore, it needs to be treated by a regularization procedure. The most general form of such a regularization procedure has been presented in \cite{CamKM12}. It allows general non-square matrices $E$ and $A$ and it can be extended to general nonlinear systems. We briefly review this regularization procedure for the linear constant coefficient case. First, we write the state equation of system \eqref{descsystem} in behavior form combining input and state to a joint vector $z=[x^T,u^T]^T$, i.\,e.,
\begin{equation}
{\mathcal E}\dot z={\mathcal A}z
\label{linbeh}
\end{equation}
with $\mathcal E=[E,\, 0]$, $\mathcal A=[A,\, B]$ partitioned accordingly.
Then following \cite{Cam87a} we form a derivative array
\begin{equation}\label{augbeh}
{\mathcal M}_\ell\dot z_\ell={\mathcal N}_\ell z_\ell,
\end{equation}
where
\[
\arraycolsep 1pt
\begin{array}{rl}
({\mathcal M}_\ell)_{i,j}=&
{i\choose j}{\mathcal E}^{(i-j)}-{i\choose j+1}{\mathcal A}^{(i-j-1)},\
i,j=0,\ldots,\ell,\\[0.12truecm]
({\mathcal N}_\ell)_{i,j}=&\left\{
\begin{array}{ll}{\mathcal A}^{(i)}~~~&\hbox{\rm for } i=0,\ldots,\ell,\ j=0,\\
                     0     &\hbox{\rm otherwise,}
\end{array}\right.\\[0.12truecm]
(z_\ell)_j=&z^{(j)},\ j=0,\ldots,\ell.
\end{array}
\]
The subsequent Theorem follows from the more general results for variable coefficient systems, see \cite{KunM06}. It connects the derivative array with the strangeness index $\mu$ and is used for index reduction.
%
%
%

\begin{theorem}\label{hyp1} Consider system (\ref{linbeh}). There exists an integer $\mu$ such that the coefficients of the derivative array \eqref{augbeh}, $(M_{\mu},N_{\mu})$, associated with $({\mathcal E},{\mathcal A})$ have the following properties, where we set
\begin{equation}
\hat{a} = a_{\mu},\quad \hat{d} = d_{\mu},\quad \hat{v} = v_0+\ldots +v_{\mu}.
\label{eq:charQuantities}
\end{equation}

\begin{enumerate}
\item  $\rank {\mathcal M}_{\mu}=(\mu+1)n-\hat{a}-\hat{v}$, i.\,e., there exists a matrix $Z$ of size
$(\mu+1)n\times (\hat{a}+\hat{v})$ and  maximal rank satisfying $Z^T{\mathcal M}_{\mu}=0$.
\item  $\rank
Z^T{\mathcal N}_{\mu}[I_{n+m}\>0\>\cdots\>0]^T=\hat{a}$, i.\,e., $Z$ can be partitioned as $Z=[\>Z_2\>\>Z_3\>]$, with $Z_2$ of
size $(\mu+1)n\times {\hat{a}}$ and $Z_3$ of size $(\mu+1)n\times {\hat{v}}$, such that
$\hat A_2=Z_2^T{\mathcal N}_{ \mu}[I_{n+m}\>0\>\cdots\>0]^T$ has full row rank
$\hat{a}$ and $Z_3^T{\mathcal N}_{ \mu}[I_{n+m}\>0\>\cdots\>0]^T=0$.
Furthermore, there exists a matrix  $T_2$ of size $\left(n+m\right)\times\left(n+m-\hat{a}\right)$ and maximal rank satisfying
$\hat A_2T_2=0$.
\item  $\rank{\mathcal E}(t)T_2 = \hat{d} = n-\hat{a}-v_{\mu}$, i.\,e., there exists a matrix $Z_1$ of
size $n\times \hat{d}$ and  maximal rank satisfying $\rank\left(\hat E_1T_2\right)=\hat{d}$ with $\hat E_1=Z_1^T{\mathcal E}$.
\end{enumerate}
Furthermore, system (\ref{linbeh}) has the same solution set as the system
\begin{equation}\label{concf}
\left [\begin{array}{c} \hat E_1 \\ 0 \\ 0 \end{array} \right ] \dot z
=\left [\begin{array}{c}
\hat A_1 \\ \hat A_2 \\ 0 \end{array} \right ]  z,
\end{equation}
where ${\hat E}_1=Z_1^T {\mathcal E}$, $\hat A_1=Z_1^T{\cal A}$ and $\hat A_2=Z_2^T{\mathcal N}_{ \mu}[I_{n+m}\>0\>\cdots\>0]^T$.
\end{theorem}

The smallest number $\mu$ for which Theorem~\ref{hyp1} holds is called the \emph{strangeness index}. The differential-algebraic system \eqref{concf} is strangeness-free, i.\,e., its strangeness index is zero. Its coefficients can be computed by using three nullspace computations, which are carried out via SVDs or $QR$ decompositions with column pivoting (cf. \cite{GolV13}) as long as this is feasible
in the available computing environment. The  system (\ref{concf}) is
a \emph{reformulation} of (\ref{linbeh}) (using the original model and its derivatives)
without changing the solution set, since no transformation of the vector $z$  has been made. The constructed submatrices $\hat A_1$ and $\hat A_2$ have been obtained
from the block matrix
\[
\left[\begin{array}{cc}
 A & B \\
 \dot A & \dot B \\
 \vdots &\vdots \\
 A^{(\mu)} & B^{(\mu)}
\end{array}\right]
\]
by transformations from the left. This has two immediate consequences \cite{KunMR01}.
First, derivatives of the input function~$u$
are nowhere needed, i.\,e., although formally
the derivatives of~$u$ occur in the derivative array, they do not
occur in the form (\ref{concf}), and hence, we do not have
any additional smoothness requirements for the input function~$u$.

Second, it follows from the construction
of $\hat A_1$ and $\hat A_2$ that the partitioning into the part
stemming from the original states~$x$ and the original controls~$u$
is not mixed up. Including the output equation,
we obtain a reformulated system of the form
\begin{subequations}\label{strfree-reform}
\begin{eqnarray}
 E_1\dot x&=&A_1x+B_1u,\label{strfree-reforma}\\
  0& = &A_2x+B_2u, \label{strfree-reformb}\\
 0& = &0,\label{strfree-reformc}\\
 y&  = & Cx+D u,\label{strfree-reformd}
\end{eqnarray}
\end{subequations}
where
\[
E_1=\hat E_1\left[\begin{array}{c} I_n \\ 0 \end{array}\right],\quad
A_i=\hat A_i\left[\begin{array}{c} I_n \\ 0 \end{array}\right],\quad
B_i=\hat A_i\left[\begin{array}{c} 0 \\ I_m \end{array}\right ],\quad i=1,2.
\]
Here $E_1,A_1$ have size $d\times n$, while $E_2,A_2$ are of size $a\times n$.
The equations in (\ref{strfree-reformc}) can just be removed from the system
and we continue with the modified model of $d+a$ equations
%
\begin{equation*}
\left [\begin{array}{c} \hat E_1 \\ 0  \end{array} \right ] \dot z
=\left [\begin{array}{c}
\hat A_1 \\ \hat A_2  \end{array} \right ]  z
\end{equation*}
together with given initial conditions.
Consistency of initial values can easily be checked, they have to satisfy the equation
\begin{equation*}
A_2 x(t_0)+B_2 u(t_0)=0,
\end{equation*}
which (if $B_2$ does not vanish) represents a restriction on the initial value of the control $u$.

In (\ref{strfree-reforma}) and (\ref{strfree-reformb}), we have
$d+a$ equations and $n$ variables in $x$ and $m$ variables in $u$.
In order for this system to be \emph{regular}, i.\,e., uniquely solvable for all sufficiently smooth inputs $u$, and all consistent initial conditions,  we would need that $d+a=n$.

If $d+a<n$, then for given $u$
we cannot expect a unique solution, i.\,e., the system is not regular and we can just attach $n-\left(d+a\right)$ variables from $x$ to $u$ and if $d+a>n$, then  we just attach $d+a-n$ of the input variables in $u$ to the vector $x$. There is freedom in the choice of the variables
that are chosen for reinterpretation, and ideally the selection should be done in such a way that the resulting descriptor system is regular if the input $u=0$ is used, but this is not necessary. Note that we must
 also change the output equation by moving appropriate  columns from $D$
to $C$ or vice versa. As a result of the reinterpretation,  we obtain a new system
\begin{eqnarray*}
 {\tilde E}_1\dot {\tilde x}&= & {\tilde A}_1\tilde x+{\tilde B}_1\tilde
 u,\\
 0 &=& {\tilde A}_2\tilde x+{\tilde B}_2\tilde
 u,\\
   y&  =&\tilde C\tilde x+\tilde D {\tilde u},
\end{eqnarray*}
where now the matrices
$\mat{c} {\tilde E}_1 \\ 0 \rix$ and $\mat{c} {\tilde A}_1 \\  {\tilde A}_2 \rix$ are square of size $\tilde n=d+a$, and
$\mat{c} {\tilde B}_1 \\ {\tilde B}_2 \rix$ is of size ${\tilde n} \times {\tilde m}$ with $\tilde m=n+m-\tilde n$.

It is often also useful to remove the feed-through term $\tilde D\tilde u$ in the output equation. This can be done by expanding the state dimension by introducing $\tilde x_{aux}\vcentcolon=\tilde D\tilde u$ and rewriting the system as
\begin{eqnarray*}
 \bar E_1 \dot {\bar x} & =& \bar A_1 {\bar x}+ \bar B_1 \tilde u,\\
 0 & =& \bar A_2 {\bar x}+ \bar B_2 \tilde u,\\
 y&=&\bar C \bar x ,
\end{eqnarray*}
with
\begin{align*}
 \bar x&=\mat{c} \tilde x \\ \tilde x_{aux} \rix,
  \bar E_1=\mat{cc} \tilde E_1 & 0 \rix,\ \bar A_1= \mat{cc} \tilde A_1 & 0 \rix,\ \bar A_2= \mat{cc}\tilde A_2 & 0 \\ 0 & I_p \rix, \
\bar B_1 = \tilde B_1,\\
\bar B_2 &= \mat{c} \tilde B_2 \\ -\tilde D \rix,\
\bar C= \mat{cc} \tilde C & I_p\rix.
\end{align*}
This method of removing the feed-through term leads to an increase of the state dimension by $p$, i.\,e., from $\tilde n$ to $\bar n=\tilde n+p$. The resulting system may again be of index higher than one as a free system with $\tilde u=0$. But in this case, see \cite{KunMR01}, there exists a linear feedback $\tilde u=K \bar x+ w$,
with $ K\in {\mathbb R}^{\tilde m,\bar n}$ such that in the closed loop system
\begin{subequations}
\begin{eqnarray}
\label{clredsysa}
\bar E \dot {\bar x}&=& (\bar A +\bar B K)\bar x+ \bar B w,\quad \bar x(t_0)=\bar x_0,\\
y&=&\bar C \bar x ,\label{clredsysb}
\end{eqnarray}
\label{clredsys}
\end{subequations}
the matrix function $(\bar A_2 +\bar B_2 K) {\bar T_2}'$ is nonsingular, and ${\bar T}_2'$ is a matrix valued function that spans the kernel of ${\bar E}_1$.
This implies that the differential-algebraic equation system in \eqref{clredsysa} is regular and of index at most one as a free system with  $w=0$, see Lemma~\ref{r-index-lemma}.
We summarize the whole regularization procedure  in the following diagram, see \cite{CamKM12}.
\vskip .5truecm
\unitlength 1.0cm
\begin{picture}(12.0,12.8)(0.0,0.0)
\arraycolsep 2pt

\put(6.0,12.8){\makebox(0,0)[ct]{
\boxed{
\begin{array}{rcl}
E \dot x&=&A x+B u,\ x(t_0)=x_0, \\
y&=& C x+D u
\end{array}
}
}}

\put(6.0,11.6){\vector(0,-1){0.5}}
\put(6.0,11.4){\makebox(0,0)[rc]{$\mu\neq 0$\quad}}
\put(6.0,11.4){\makebox(0,0)[lc]{\quad index reduction in behavior}}

\put(6.0,11.0){\makebox(0,0)[ct]{
\boxed{
\begin{array}{rcl}
E_1\dot x&=&A_1x+B_1u,\\
0&=&A_2x+B_2u,\\
0&=&0,\\
y&=&Cx+D u
\end{array}
}
}}

\put(6.0,8.8){\vector(0,-1){0.5}}
\put(6.0,8.6){\makebox(0,0)[rc]{\quad}}
\put(6.0,8.6){\makebox(0,0)[lc]{\quad remove $0=0$ eq.}}

\put(6.0,8.2){\vector(0,-1){0.5}}
\put(6.0,8.0){\makebox(0,0)[rc]{$0 \neq  A_2 x_0 +B_2 u(t_0)$\quad}}
\put(6.0,8.0){\makebox(0,0)[lc]{\quad  cond. for consistency}}

\put(6.0,7.6){\vector(0,-1){0.5}}
\put(6.0,7.4){\makebox(0,0)[rc]{$a+d\neq n$\quad}}
\put(6.0,7.4){\makebox(0,0)[lc]{\quad reinterpret variables}}

\put(6.0,7.0){\makebox(0,0)[ct]{
\boxed{
\begin{array}{rcl}
{\tilde E}_1\dot {\tilde x}&=&{\tilde A}_1\tilde x+{\tilde B}_1\tilde u,\\
0&=&{\tilde A}_2\tilde x+{\tilde B}_2\tilde u,\\
\tilde y&=&\tilde C\tilde x+\tilde D {\tilde u}
\end{array}
}
}}

\put(6.0,5.2){\vector(0,-1){0.5}}
\put(6.0,5.0){\makebox(0,0)[rc]{$\tilde D \neq 0$\quad}}
\put(6.0,5.0){\makebox(0,0)[lc]{\quad remove feed-through}}

\put(6.0,4.7){\makebox(0,0)[ct]{
\boxed{
\begin{array}{rcl}
{\bar E}_1\dot {\bar x}&=&{\bar A}_1\bar x+{\bar B}_1\bar u,\\
0&=&{\bar A}_2\bar x+{\bar B}_2\bar u,\\
y&=&\bar C\bar x
\end{array}
}
}}

\put(6.0,3.0){\vector(0,-1){0.5}}
\put(6.0,2.8){\makebox(0,0)[rc]{not strangen.-free for $u=0$\quad}}
\put(6.0,2.8){\makebox(0,0)[lc]{\quad perform feedback $\bar u=K \bar x+w$}}

\put(6.0,2.5){\makebox(0,0)[ct]{
\boxed{
\begin{array}{rcl}
{\bar E}_1\dot {\bar x}&=&({\bar A}_1+{\bar B}_1K)\bar x+{\bar B}_1w,\\
0&=&({\bar A}_2+{\bar B}_2K)\bar x+{\bar B}_2w,\\
\bar y&=&\bar C\bar x.
\end{array}
}
}}

\end{picture}

\noindent In the following we assume that the system has been regularized to the form \eqref{clredsys} and furthermore that $\bar B$ and $\bar C^T$ have full column rank.  Otherwise, we can just reduce the input vector
or the output vector. In abuse of notation, we denote the resulting system again
in the original notation
\begin{equation}
\begin{aligned}
E \dot {x}&= A x+ B u,\quad x(t_0)= x_0,\\
 y&=C x ,
 \end{aligned}
 \label{eq:reg}
\end{equation}
Note that the resulting system may not satisfy the desired controllability and observability  conditions associated with the finite spectrum, and even if, then it may be close to a system that does not satisfy these conditions. To remove uncontrollable and unobservable finite parts, i.\,e., to make the system minimal,
some further transformations may be necessary. In the following section we discuss condensed forms under
orthogonal transformations which can be used to check all the controllability conditions from Section~\ref{controllability}.

\section{Condensed Forms}\label{condensed}

To verify the controllability and observability conditions,
using equivalence transformations such as (\ref{a-25}), the regularized system \eqref{eq:reg} is transformed to a condensed form that reveals these properties. The following condensed form has been presented in full generality in \cite{BunBMN99}. It uses only real orthogonal
transformations and can be computed using
algorithms that are numerically stable in the sense that
in finite precision arithmetic, the computed condensed
form is what would have been obtained using exact arithmetic
from a rounding-error-small perturbation of the original
descriptor system. In the following  we adopt the notation
that a matrix $\Sigma_j$ is a non-singular $j\mbox{-by-}j$ diagonal matrix,
and $0$ denotes the null-matrix of any size.

Unfortunately, all condensed forms rely on numerical rank decisions
of transformed submatrices of \(E\), \(A\), \(B\) and \(C\).
This is a serious problem, since arbitrarily small perturbations
of a rank deficient matrix may change its rank.

\begin{theorem}\label{thm3.1}\cite{BunMN94}
Let  $E,\,A\in \mathbb R^{n,n}$, $B\in \mathbb R^{n,m}$, and $C \in \mathbb R^{p,n}$,
  where $B$ and $C$ are of full column and row rank, respectively.
  Then, there exist orthogonal matrices $U,\,V \in \mathbb R^{n,n}$,
  $W\in \mathbb R^{m,m}$, and $Y \in \mathbb R^{p,p}$ such that
  \begin{subequations}
  \begin{eqnarray}
       \label{3.2a}
       U^T E V & = & \bmat{
                            &    t_1   & n-t_1 \cr
                        t_1 & \Sigma_{t_1} &   0   \cr
                      n-t_1 &    0     &   0   \cr
                 }
                 ,\\
       \label{3.2b}
        U^T B W & = &\bmat{
                        &      k_1    & k_2     \cr
                    t_1 &      B_{11} & B_{12}  \cr
                    t_2 &      B_{21} & 0  \cr
                    t_3 &      B_{31} & 0  \cr
                  n-t_1-t_2-t_3 & 0    &  0  \cr
                  }
                  ,
        \\
        \label{3.2c}
        Y^T C V &=& \bmat{
                         &   t_1  &  s_2 & t_5 & n-t_1-s_2-t_5       \cr
                  \ell_1 & C_{11} & {C}_{12} & C_{13} & 0\cr
                  \ell_2 & C_{21} &      0  & 0 & 0
                  }
                  ,
        \\
        \label{3.2d}
        U^T A V & = &
        \bmat{
                & t_1 & s_2    & t_5    & t_4    & t_3    & s_6    \cr
 t_1 & A_{11} & A_{12} & A_{13} & A_{14} & A_{15} & A_{16} \cr
 t_2 & A_{21} & A_{22} &      A_{23} &      A_{24} &  0     &  0  \cr
 t_3 &A_{31} & A_{32} & A_{33} &      A_{34} & \Sigma_{t_3} &  0  \cr
 t_4 &A_{41} & A_{42} &  A_{43} & \Sigma_{t_4} &       0     &  0  \cr
 t_5  &A_{51} &  0 & \Sigma_{t_5} &       0     &       0     &  0  \cr
 t_6 &A_{61} &  0   &   0     &       0     &       0     &  0  \cr
        }
        .
\end{eqnarray}
\end{subequations}
The matrix \(B_{12}\) has full column rank,
\(C_{21}\) has full row rank, and the matrices
\begin{eqnarray*}
   \mat {c} B_{21} \\* B_{31} \rix \in \C^{k_1 \times k_1},
   & \quad &
   \mat{cc} C_{12}  &  C_{13} \rix \in \C^{\ell_1 \times \ell_1}
\end{eqnarray*}
are square and non-singular and are of dimension $k_1=t_2+t_3$
and $\ell_1=s_2+t_5$, respectively. Here
\(t_j\), \(s_j\), \(k_j\) and \(\ell_j\) are non-negative integers
displaying the number of rows or columns in the corresponding block row or
column of the matrices.
A zero value of one of these integers indicates that the corresponding
block row or column does not appear.

\end{theorem}

\noindent As a corollary we can characterize controllability and observability conditions of Section~\ref{controllability}.
\begin{corollary}\label{cor3.9}
Consider a  system of the form \eqref{eq:reg} and let the system be transformed to the condensed form (\ref{3.2a})--(\ref{3.2d})
of Theorem~\ref{thm3.1}.
\begin{enumerate}
\item
        The pair $(E,A)$ is regular and of index at most one if and only
        if $s_6=t_6=0$ and $A_{22}$
        is nonsingular.
\item
        Condition {\bf C2} holds
        if and only if $t_6=0$.
\item
        Condition {\bf O2} holds
        if and only if $s_6=0$.
\item
        $\rank [ E,B]=t_1+t_2+t_3$, and thus
        $\rank [ E,B]=n$ if and only if $t_4=t_5=t_6=0$.
\item
        $\rank \mat{c} E \\ C \rix=t_1+s_2+t_5$, and thus
        $\rank \mat{c} E \\ C \rix=n$ if and only if $t_4=t_3=s_6=0$.
\item
        $
                \rank \mat{cc} E & B \\ C & 0 \rix
                =
                t_1+t_2+s_2+t_3+t_5+\min(\ell_2,k_2)
        $.
\end{enumerate}
\end{corollary}

\noindent If we have computed the condensed form from the regularized system \eqref{eq:reg} then we should have $t_6=s_6=0$ and the system is of index at most one as a free system. The staircase form allows to check whether the regularization procedure has been successful.

\section{{\tt Matlab} Functions in detail}\label{sec:Matlab}

In the following several {\tt Matlab} functions are presented, which create a regularized Loewner realization based on tangential interpolation data of the transfer function. It strongly builds on the procedure of \cite{MayA07}, see Section \ref{sec:realization}, followed by a regularization based on the methods and results outlined in the previous sections. As a result, we obtain a realization which is regular and strangeness-free as well as completely controllable and observable.

\subsection{Realization}\label{sec:realizationFunc}

\subsubsection*{Syntax}

\verb+[E,A,B,C,D,mu,la,V,W,L,R] = realization(S,F)+

\noindent\verb+[E,A,B,C,D,mu,la,V,W,L,R,U_trans,V_trans,W_trans,Y_trans,...+

\verb+ L_trans,R_trans,Feedb] = realization(S,F)+

\noindent\verb+[E,A,B,C,D,mu,la,V,W,L,R,U_trans,V_trans,W_trans,Y_trans,...+

\verb+ L_trans,R_trans,Feedb] = realization(S,F,tol)+

\noindent\verb+[E,A,B,C,D,mu,la,V,W,L,R,U_trans,V_trans,W_trans,Y_trans,...+

\verb+ L_trans,R_trans,Feedb] = realization(S,F,tol,sindexflag)+

\subsubsection*{Arguments}

\itxt{function {\tt realization}}
\bal
\verb+S+          &	Vector of length $N$ of interpolation points (is split into two disjoint \\
				  &	interpolation point sets \verb+mu+ and \verb+lambda+)\\
\verb+F+          &	Array which contains the transfer function values at points \verb+S+; this is either\\
				  & a vector of length N (scalar interpolation case) or a $p\times m\times N$ array\\
				  & consisting of $N$ $p\times m$ matrices (matrix interpolation case)\\
\verb+tol+        & scalar specifying the tolerance value for rank decisions in function \verb+hypo+\\
				  & (default: \verb+10*eps+, where \verb+eps+ is the floating-point relative accuracy $2^{-52}$)\\
\verb+sindexflag+ & boolean (default: \verb+true+); if \verb+true+, index reduction and regularization are\\
				  & performed; if \verb+false+, Loewner realization is provided without\\
				  & post-processing steps
\eal

\otxt{function {\tt realization}}
\bal
\verb+[E,A,B,C,D]+				 &	matrices corresponding to a system which is regularized and\\
							     &  strangeness-free (if \verb+sindexflag+ is set to \verb+true+) and whose transfer\\
								 &	function $H(z)=C(zE-A)^{-1}B+D$ interpolates the data (\verb+S+,\verb+F+)\\
\verb+mu+						 &	vector of size $\nu$ containing the left interpolation points $\mu_j$ (cf. \eqref{leftintcond})\\
\verb+la+						 &	vector of size $\rho$ containing the right interpolation points $\lambda_i$ (cf. \eqref{rightintcond})\\
\verb+V+						 &	\textit{scalar interpolation case}: vector of size $\nu$ containing the left\\
 								 &	interpolation values $v_j$ belonging to $\mu_j$  (cf.  \eqref{leftintcond} with $l_j=1$)\\
								 &	\textit{matrix interpolation case}: matrix of dimension $\nu\times m$ containing\\
								 &  the left interpolation values $v_j\in\C^{1,m}$ (as rows) generated by\\
								 &  random left tangential directions $l_j\in\C^{1,p},\, j=1,\ldots,\nu$ (cf. \eqref{leftintcond})\\
\verb+W+						 &	\textit{scalar interpolation case}: vector of size $\rho$ containing the right\\
 								 &	interpolation values $w_i$ belonging to $\lambda_i$  (cf.  \eqref{rightintcond} with $r_i=1$)\\
								 &	\textit{matrix interpolation case}: matrix of dimension $p\times\rho$ containing the\\
								 &  right interpolation values $w_i\in\C^{p}$ (as columns) generated by\\
								 &  random right tangential directions $r_i\in\C^{m},\, i=1,\ldots,\rho$ (cf. \eqref{rightintcond})\\
\verb+L+						 &	matrix of dimension $\nu\times p$ containing (random) left tangential\\
								 &  directions $l_j\in\C^{1,p}$ as rows (set to one in scalar interpolation case)\\
\verb+R+						 &	matrix of dimension $m\times \rho$ containing (random) right tangential\\
								 &  directions $r_i\in\C^{m}$ as columns (set to one in scalar interpolation\\
								 &  case)\\
\verb+U_trans+, \verb+V_trans+,  &	\\
\verb+W_trans+, \verb+Y_trans+	 &  matrices corresponding to the transformation matrices $U$, $V$, $W$,\\
 								 &  and $Y$ of Theorem~\ref{thm3.1} (see also description of function \verb+staircase+)\\
\verb+L_trans+, \verb+R_trans+	 &	transformation matrices from block Gaussian elimination in\\
								 &  function \verb+regularization+\\
\verb+Feedb+					 &	feedback matrix from function \verb+regularization+ to make\\ 						
								 &	$(E,A)$ regular
\eal

\begin{remark}{\rm
The vectors $la$ and $mu$ together form $S$ such that $\nu+\rho=N$ and the sizes of $la$ and $mu$ differ by one when $N$ is odd and are the same when $N$ is even. Furthermore, if $S$ contains values with non-zero imaginary part the complex conjugate values are added to $S$ to ensure that the output realization \verb+[E,A,B,C,D]+ consists of real-valued matrices.}
\end{remark}

\subsubsection*{Description}

\verb+[E,A,B,C,D,mu,la,V,W,L,R] = realization(S,F)+ constructs matrices $E$, $A$, $B$, $C$, and $D$ such that the transfer function $H(s)=C(sE-A)^{-1}B+D$ interpolates the given data as described in Section \ref{sec:realization}. First, the Loewner matrix $\LL$ and the shifted Loewner matrix $\LL_\sigma$ as well as the corresponding matrices $V,W,L,R$ and vectors $\mu,\lambda$ are constructed by use of the {\tt{Matlab}} function
\noindent\verb+[LL,sLL,mu,la,V,W,L,R] = loewner_mod(S,F)+. Then two cases have to be considered. If $\rho=\nu$ and $\det(\tilde s\LL-\LL_\sigma)\neq0$ for all $\tilde s\in\{\lambda_i\}\cup\{\mu_j\}$, which implies that $\tilde s\LL-\LL_\sigma$ is quadratic and nonsingular (regular case), the Loewner realization $\left(E,A,B,C\right)$ is given by $E=-\LL$, $A=-\LL_\sigma$, $B=V$, and $C=W$ resulting in a system with the desired interpolation properties. The second case is the nonregular case. To make sure that we still get a regular system, we need to ensure that
\begin{equation}
\rank(\tilde s\LL-\LL_\sigma)=\rank\mat{cc}\LL & \LL_\sigma\rix=\rank\mat{c}\LL \\* \LL_\sigma\rix,\quad \mbox{\rm for all } \tilde s\in\{\lambda_i\}\cup\{\mu_j\}.
\label{eq:rankCondition}
\end{equation}
If this condition is satisfied, we choose an arbitrary $\tilde s\in\{\lambda_i\}\cup\{\mu_j\}$ and compute the skinny SVD $(\tilde s\LL-\LL_\sigma) = Y\Sigma X^T$, see \cite{GolV13}  with a nonsingular diagonal matrix $\Sigma$ and transformation matrices $Y$ and $X$ with pairwise orthonormal columns. In this case the Loewner realization, see \cite{MayA07}, is given by
\[
E=-Y^T\LL X,\quad A=-Y^T\LL_\sigma X,\quad B=Y^TV,\quad C=WX,\quad D=0.
\]

Even if the regularity of the matrix pencil $\left(E,A\right)$ is guaranteed analytically, in the finite precision case, we cannot be sure about this. Further important properties as the index, controllability, and observability are also unknown, in general. Thus, to obtain a regular, strangeness-free, completely observable, and completely controllable system, some further steps have to be performed (only executed if \verb+sindexflag+ is true).

First, the index is reduced by applying Theorem \ref{hyp1} using the function \noindent\verb+hypo+. The resulting system is strangeness-free and is of the form \eqref{strfree-reform}. The vanishing equations can be neglected such that the number of equations decreases to $a+d$. If this number differs from the number of variables, either some of the components of $x$ have to be attached to the vector $u$ or vice versa. This changes the input dimension $m$ such that the size of the transfer function does not fit to the tangential interpolation data anymore. However, if the problem is well-posed, this case should not occur.

To obtain more insight into the controllability and observability properties of the realization, the function \verb+staircase+ is called, which computes the condensed form of the realization $\left(E,A,B,C\right)$ according to Theorem \ref{thm3.1}. The system matrices corresponding to the condensed form are denoted with $EC$, $AC$, $BC$, and $CC$, respectively. If $B$ does not have full column rank or if $C$ does not have full row rank, the resulting zero columns in $BC$ or zero rows in $CC$ are canceled which decreases the dimensions of $L$, $R$, $V$, and $W$ accordingly.

The subsequent regularization procedure, performed by the function \verb+regularization+, eliminates the non-controllable and non-observable parts leading to system matrices $E_{reg}$, $A_{reg}$, $B_{reg}$, and $C_{reg}$ as
\begin{equation*}
E_{reg}= \begin{bmatrix}
E_{11} & \\
	& 0
\end{bmatrix},\quad
A_{reg}=\begin{bmatrix}
A_{11} & \\
	& A_{22}
\end{bmatrix},
\end{equation*}

\begin{equation*}
B_{reg}=\begin{bmatrix}
B_{11} & B_{12}\\
B_{21}	& 0
\end{bmatrix},\quad
C_{reg}=\begin{bmatrix}
C_{11} & C_{12} \\
C_{21}	& 0
\end{bmatrix},
\end{equation*}
where the blocks $E_{11}$ and $A_{22}$ are nonsingular. Consequently, the pencil $(\tilde sE_{reg}-A_{reg})$ is regular. Finally, expanding the product of the matrices, one gets

\begin{equation*}
\begin{bmatrix}
C_{11}\\
C_{21}
\end{bmatrix}\left(sE_{11}-A_{11}\right)^{-1}
\begin{bmatrix}
B_{11} & B_{12}
\end{bmatrix}
-C_{12}A_{22}^{-1}B_{21},
\end{equation*}
such that we can decrease the state-space dimension of the realization by setting
\begin{equation*}
E=E_{11}, \quad A=A_{11}, \quad
B=\begin{bmatrix}
B_{11} & B_{12}
\end{bmatrix}, \quad
C=\begin{bmatrix}
C_{11}\\
C_{21}
\end{bmatrix}, \quad
D=-C_{12}A_{22}^{-1}B_{21}.
\end{equation*}

\subsection{Hypo}

\subsubsection*{Syntax}

\verb+[E1_hat,A1_hat,A2_hat,d,a,v,mu_max,sig] = hypo(E,A,mu,tol,varargin)+

\subsubsection*{Arguments}

\itxt{function {\tt hypo}}
\bal
\verb+E+			&	matrix $\mathcal{E}$ of the system's behavior form as in \eqref{linbeh}, i.\,e., $\mathcal{E}=[E,\, 0]$\\
\verb+A+			&	matrix $\mathcal{A}$ of the system's behavior form as in \eqref{linbeh}, i.\,e., $\mathcal{A}=[A,\, B]$\\
\verb+mu+			&	corresponds to the index $\ell$ of the inflated system \eqref{augbeh} (default: $0$)\\
\verb+tol+			&   scalar specifying the tolerance value for rank decisions (default: \verb+10*eps+)\\
\verb+varargin+  	&	contains $v_0,\ldots,v_{\ell-1}$ (default: empty)
\eal

\otxt{function {\tt hypo}}
\bal
\verb+E1_hat+, \verb+A1_hat+,	&	\\
\verb+A2_hat+				    &	blocks of the reformulated system (\ref{concf})\\
\verb+d+						&	number of differential equations ($\hat{d}$)\\
\verb+a+						&	number of algebraic equations ($\hat{a}$)\\
\verb+v+						&	number of vanishing equations ($v_{\mu}$)\\
\verb+mu_max+					&	strangeness index of the original system $\left(\mathcal E, \mathcal A\right)$\\
\verb+sig+						&	error resulting from rank decision
\eal

\subsubsection*{Description}

The function \verb+hypo+ successively inflates the system $\mathcal{E}\dot z=\mathcal{A}z$ by differentiation which leads to inflated systems $\mathcal{M}_\ell \dot z_\ell=\mathcal{N}_\ell z_\ell$ with $\ell$ starting at $0$ and being incremented by one in each step. This procedure is continued until the rank conditions of Theorem \ref{hyp1} are fulfilled yielding $\mu,\,\hat{d},\,\hat{a}$, and $\hat{v}$.
The matrix $Z$ is computed by means of an SVD of $\mathcal{M}_\ell$ using those left singular vectors that lie in the left null space of $\mathcal{M}_\ell$. The first singular value that is considered to be negligibly small (during rank decision based on \verb+tol+) is used as an error measurement of the procedure.

The matrices $Z_2$ and $T_2$ are determined based on the SVD
\begin{equation*}
Z^T\mathcal{N}_\ell \mat{cccc} I & 0 & \cdots & 0\rix^T=USV^T.
\end{equation*}
$Z_2$ consists of the first $\hat{a}$ columns of $ZU$, i.\,e., such that $Z_2^T\mathcal{N}_\ell \mat{cccc} I & 0 & \cdots & 0\rix^T$ has full row rank and $T_2$ consists of those columns of $V$ lying in the right null space of the matrix $Z_2^T\mathcal{N}_\ell \mat{cccc} I & 0 & \cdots & 0\rix^T$. The difference between the number of columns of $Z$ and the number of algebraic constraints $\hat{a}$ is equal to $\hat{v}$, cf. Theorem \ref{hyp1}.

Finally, $Z_1$ is determined by calculating a QR-decomposition of $\mathcal{E}T_2$ and by choosing $\hat d$ columns of Q such that $Z_1^T \mathcal{E}T_2$ has full rank $\hat d$. If the sum $\hat{d}+\hat{a}+v_{\ell}$ (using $v_{\ell} = \hat{v}-\sum_{i=0}^{\ell-1}v_i$, cf. (\ref{eq:charQuantities})) differs from the number of equations of the system $\mathcal{E}\dot z=\mathcal{A}z$, the index $\ell$ is increased by one and \verb+hypo+ is called with \verb+varargin+ containing $v_0,\ldots,v_{\ell-1}$. Otherwise the index reduction is complete and we set $\hat E_1=Z_1^T\mathcal{E}$, $\hat A_1=Z_1^T\mathcal{A}$ and $\hat A_2=Z_2^T\mathcal{N}_\ell \mat{cccc} I & 0 & \cdots & 0\rix^T$. The number $\mu_{max}$ corresponds to the smallest index $\ell$ needed to satisfy Theorem \ref{hyp1}. This number is equal to the strangeness index of the original system $\left(\mathcal E, \mathcal A\right)$.

\subsection{Staircase}

\subsubsection*{Syntax}

\verb+[EC,AC,BC,CC,U,V,W,Y,t,s,k,l] = staircase(E,A,B,C)+

\subsubsection*{Arguments}

\itxt{function {\tt staircase}}
\bal
\verb+E+	&	$n1\times n2$ matrix\\
\verb+A+	&	$n1\times n2$ matrix\\
\verb+B+	&	$n1\times m$ matrix\\
\verb+C+	&	$p\times n2$ matrix
\eal

\otxt{function {\tt staircase}}
\bal
\verb+EC+, \verb+AC+, \verb+BC+, \verb+CC+	&	condensed form of the input system matrices $\left(E, A, B, C\right)$ according to \\
											&   Theorem \ref{thm3.1}\\
\verb+U+, \verb+V+, \verb+W+, \verb+Y+		&	orthogonal matrices that transform $E$, $A$, $B$, and $C$ to condensed form,\\ 				
											&	i.\,e., $EC=U^TEV$, $AC=U^TAV$, $BC=U^TBW$, and $CC=Y^TCV$\\
\verb+t+, \verb+s+, \verb+k+, \verb+l+		&	vectors containing the block dimensions of the condensed form, see\\
											&   Theorem \ref{thm3.1}
\eal

\subsubsection*{Description}

The algorithm follows the constructive proof of Theorem \ref{thm3.1}, which is presented in \cite{BunMN94}. For that, numerous SVDs are used to transform the input matrices $E$, $A$, $B$, and $C$ into the form

\begin{equation}
\begin{aligned}
\verb+EC+&=\begin{bmatrix}
\Sigma_{E} & 0 & 0 & 0 & 0 & 0\\
0 & 0 & 0 & 0 & 0 & 0\\
0 & 0 & 0 & 0 & 0 & 0\\
0 & 0 & 0 & 0 & 0 & 0\\
0 & 0 & 0 & 0 & 0 & 0\\
0 & 0 & 0 & 0 & 0 & 0
\end{bmatrix},\quad
&\verb+AC+=\begin{bmatrix}
A_{11} & A_{12} & A_{13} & A_{14} & A_{15} & A_{16}\\
A_{21} & A_{22} & A_{23} & A_{24} & 0 & 0\\
A_{31} & A_{32} & A_{33} & A_{34} & \Sigma_{t_3} & 0\\
A_{41} & A_{42} & A_{43} & \Sigma_{t_4} & 0 & 0\\
A_{51} & 0 & \Sigma_{t_5} & 0 & 0 & 0\\
A_{61} & 0 & 0 & 0 & 0 & 0
\end{bmatrix},\\
\verb+BC+&=\begin{bmatrix}
B_{11} & B_{12} & 0\\
B_{21} & 0 & 0\\
B_{31} & 0 & 0\\
0 & 0 & 0\\
0 & 0 & 0\\
0 & 0 & 0
\end{bmatrix},\quad
&\verb+CC+=\begin{bmatrix}
C_{11} & C_{12} & C_{13} & 0 & 0 & 0\\
C_{21} & 0 & 0 & 0 & 0 & 0\\
0 & 0 & 0 & 0 & 0 & 0
\end{bmatrix},
\end{aligned}
\label{eq:generalCondensedForm}
\end{equation}
where $\verb+EC+$ and $\verb+AC+$ are of size $(t_1+t_2+t_3+t_4+t_5+t_6)\times (t_1+s_2+t_5+t_4+t_3+s_6)$,
$\verb+BC+$ is of size $(t_1+t_2+t_3+t_4+t_5+t_6)\times (k_1+k_2+(m-k_1-k_2))$ and $\verb+CC+$ is of size $(l_1+l_2+(p-l_1-l_2))\times (t_1+s_2+t_5+t_4+t_3+s_6)$. Accordingly, we have $n1=t_1+t_2+t_3+t_4+t_5+t_6$ and $n2=t_1+s_2+t_5+t_4+t_3+s_6$. Note that the difference between \eqref{eq:generalCondensedForm} and the condensed form presented in Theorem \ref{thm3.1} is that \eqref{eq:generalCondensedForm} allows for general input matrices $B$ and $C$ without assuming full row or column rank. During the algorithm also the transformation matrices are built such that

\begin{equation*}
\verb+EC+=U^TEV,\quad \verb+AC+=U^TAV, \quad \verb+BC+=U^TBW,\quad \verb+CC+=Y^TCV,
\end{equation*}
with $U\in\C^{n1\times n1}$, $V\in\C^{n2\times n2}$, $W\in\C^{m\times m}$, and $Y\in\C^{p\times p}$.

\subsection{Regularization}

\subsubsection*{Syntax}

\verb+[E,A,B,C,L_trans,R_trans,Feedb]= regularization(EC,AC,BC,CC,t,s,k,l)+

\subsubsection*{Arguments}

\itxt{function {\tt regularization}}
\bal
\verb+EC+, \verb+AC+, \verb+BC+, \verb+CC+	&	matrices in condensed form generated by the function \verb+staircase+\\
\verb+t+, \verb+s+, \verb+k+, \verb+l+		&	vectors containing the block dimensions of the condensed form generated\\
											&   by the function \verb+staircase+
\eal

\otxt{function {\tt regularization}}
\bal
\verb+E+, \verb+A+, \verb+B+, \verb+C+		&	controllable and observable system where $(E,A)$ is regular\\
\verb+L_trans+, \verb+R_trans+				&	left and right transformation matrices such that\\						
											&	$E=L_{trans}\verb+EC+R_{trans}$, $A=L_{trans}\verb+AC+R_{trans}+L_{trans}\verb+BC+Feedb$,\\
											&   $B=L_{trans}\verb+BC+$ and $C=\verb+CC+R_{trans}$\\
\verb+Feedb+								&	feedback matrix, which ensures that the block $A_{22}$ of $A$ is\\
										    &   nonsingular
\eal

\subsubsection*{Description}

In the function \verb+regularization+, first it is checked whether the input system can be made regular and of index one. This means that we have to ensure, that the matrices $\verb+EC+$ and $\verb+AC+$ are quadratic and that $t_6=s_6=0$. If this is true, the matrices have the following form:

\begin{equation*}
\verb+EC+=\begin{bmatrix}
\Sigma_{E} & 0 & 0 & 0 & 0\\
0 & 0 & 0 & 0 & 0 \\
0 & 0 & 0 & 0 & 0 \\
0 & 0 & 0 & 0 & 0 \\
0 & 0 & 0 & 0 & 0
\end{bmatrix},\quad
\verb+AC+=\begin{bmatrix}
A_{11} & A_{12} & A_{13} & A_{14} & A_{15}\\
A_{21} & A_{22} & A_{23} & A_{24} & 0\\
A_{31} & A_{32} & A_{33} & A_{34} & \Sigma_{t_3}\\
A_{41} & A_{42} & A_{43} & \Sigma_{t_4} & 0\\
A_{51} & 0 & \Sigma_{t_5} & 0 & 0
\end{bmatrix},
\end{equation*}

\begin{equation*}
\verb+BC+\begin{bmatrix}
B_{11} & B_{12}\\
B_{21} & 0\\
B_{31} & 0\\
0 & 0\\
0 & 0
\end{bmatrix},\quad
\verb+CC+=\begin{bmatrix}
C_{11} & C_{12} & C_{13} & 0 & 0\\
C_{21} & 0 & 0 & 0 & 0
\end{bmatrix},
\end{equation*}
where the block $A_{22}$ is quadratic ($s_2\overset{!}{=}t_2$) and the zero columns and rows of $\verb+BC+$ and $\verb+CC+$ are canceled out in the function \verb+realization+.

The blocks $\Sigma_{t_i}$ in $\verb+AC+$ are invertible diagonal matrices such that a block Gaussian elimination can be performed to eliminate the corresponding rows and columns inside $\verb+AC+$ leading to
\begin{equation*}
\verb+A1+=\begin{bmatrix}
A_{11} & A_{12} & 0 & 0 & 0\\
A_{21} & A_{22} & 0 & 0 & 0\\
0 & 0 & 0 & 0 & \Sigma_{t_3}\\
0 & 0 & 0 & \Sigma_{t_4} & 0\\
0 & 0 & \Sigma_{t_5} & 0 & 0
\end{bmatrix}.
\end{equation*}
$\verb+BC+$ and $\verb+CC+$ are transformed accordingly to $\verb+B1+$ and $\verb+C1+$ without changing the block structure while $\verb+EC+$ stays completely unchanged due to its zero-block structure. In the end we are only interested in the system's transfer function $H(s)=C(sE-A)^{-1}B+D$. Thus, we can restrict ourselves to the upper left $2\times 2$ block of $\verb+EC+$ and $\verb+A1+$, since by multiplying the lower right part of $(sE-A)^{-1}$, namely the block
\[
\begin{bmatrix}
0& 0& -\Sigma_{t_5}^{-1}\\
0 & -\Sigma_{t_4}^{-1} & 0\\
-\Sigma_{t_3}^{-1} &0 &0
\end{bmatrix},
\]
by the corresponding blocks of $\verb+B1+$ and $\verb+C1+$, it cancels out and, hence, it does not contribute to the transfer function. The system can be reduced to

\begin{equation}
\begin{aligned}
E_{new}&=\begin{bmatrix}
E_{11} & 0\\
0 & 0
\end{bmatrix}, \quad
&A_{new}=\begin{bmatrix}
A_{11} & A_{12}\\
A_{21} & A_{22}
\end{bmatrix},\\
B_{new}&=\begin{bmatrix}
B_{11} & B_{12}\\
B_{21} & 0
\end{bmatrix},\quad
&C_{new}=\begin{bmatrix}
C_{11} & C_{12}\\
C_{21} & 0
\end{bmatrix},
\end{aligned}
\label{eq:2by2System}
\end{equation}
where, by abuse of notation, we have redefined the naming of the matrix blocks, i.\,e., $A_{11}$ in \eqref{eq:2by2System} is not necessarily the same as $A_{11}$ in \eqref{eq:generalCondensedForm} and so on.

If the block $A_{22}$ is singular, then the pencil $(sE_{new}-A_{new})$ will not be strangeness-free. In this case a feedback is added using the fact that the block $B_{21}$ is invertible by its construction in \verb+staircase+. We construct a matrix $Feedb$ such that
\[
A_{new}+B_{new}Feedb=\begin{bmatrix}
\tilde{A}_{11} & \tilde{A}_{12}\\
0 & \sigma I
\end{bmatrix}=A_{new2},
\]
where $I$ is the identity matrix and $\sigma$ denotes the smallest singular value of the first block row of $A_{new}$.
Using block Gaussian elimination we can then  transform $A_{new2}$ into block diagonal form and obtain the desired regularized system together with the transformation matrices $L_{trans}$ and $R_{trans}$.

\section{Numerical Example \label{sec:num}}

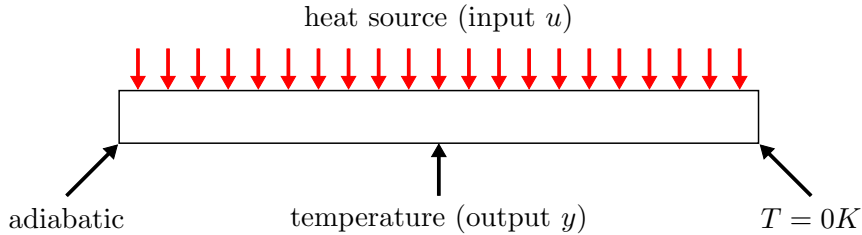
\begin{figure}
\centering
\begin{tikzpicture}
	\coordinate (a) at (0,0);
	\draw[line width = 0.5mm, color = red, preaction={-triangle 60,thin,draw,shorten >=-0.5mm,color=red}] ($(a)+(-4,0)$) -- ($(a)+(-4,-0.5)$);
	\draw[line width = 0.5mm, color = red, preaction={-triangle 60,thin,draw,shorten >=-0.5mm,color=red}] ($(a)+(-3.6,0)$) -- ($(a)+(-3.6,-0.5)$);
	\draw[line width = 0.5mm, color = red, preaction={-triangle 60,thin,draw,shorten >=-0.5mm,color=red}] ($(a)+(-3.2,0)$) -- ($(a)+(-3.2,-0.5)$);
	\draw[line width = 0.5mm, color = red, preaction={-triangle 60,thin,draw,shorten >=-0.5mm,color=red}] ($(a)+(-2.8,0)$) -- ($(a)+(-2.8,-0.5)$);
	\draw[line width = 0.5mm, color = red, preaction={-triangle 60,thin,draw,shorten >=-0.5mm,color=red}] ($(a)+(-2.4,0)$) -- ($(a)+(-2.4,-0.5)$);
	\draw[line width = 0.5mm, color = red, preaction={-triangle 60,thin,draw,shorten >=-0.5mm,color=red}] ($(a)+(-2,0)$) -- ($(a)+(-2,-0.5)$);
	\draw[line width = 0.5mm, color = red, preaction={-triangle 60,thin,draw,shorten >=-0.5mm,color=red}] ($(a)+(-1.6,0)$) -- ($(a)+(-1.6,-0.5)$);
	\draw[line width = 0.5mm, color = red, preaction={-triangle 60,thin,draw,shorten >=-0.5mm,color=red}] ($(a)+(-1.2,0)$) -- ($(a)+(-1.2,-0.5)$);
	\draw[line width = 0.5mm, color = red, preaction={-triangle 60,thin,draw,shorten >=-0.5mm,color=red}] ($(a)+(-0.8,0)$) -- ($(a)+(-0.8,-0.5)$);
	\draw[line width = 0.5mm, color = red, preaction={-triangle 60,thin,draw,shorten >=-0.5mm,color=red}] ($(a)+(-0.4,0)$) -- ($(a)+(-0.4,-0.5)$);
	\draw[line width = 0.5mm, color = red, preaction={-triangle 60,thin,draw,shorten >=-0.5mm,color=red}] (a) -- ($(a)+(0,-0.5)$) node[midway,above=0.3cm,color=black] {heat source (input $u$)} node[midway,below=0.3cm,draw,line width = 0.2mm,color=black,minimum width=8.5cm,minimum height=0.7cm] (b) {};
	\draw[line width = 0.5mm, color = red, preaction={-triangle 60,thin,draw,shorten >=-0.5mm,color=red}] ($(a)+(0.4,0)$) -- ($(a)+(0.4,-0.5)$);
	\draw[line width = 0.5mm, color = red, preaction={-triangle 60,thin,draw,shorten >=-0.5mm,color=red}] ($(a)+(0.8,0)$) -- ($(a)+(0.8,-0.5)$);
	\draw[line width = 0.5mm, color = red, preaction={-triangle 60,thin,draw,shorten >=-0.5mm,color=red}] ($(a)+(1.2,0)$) -- ($(a)+(1.2,-0.5)$);
	\draw[line width = 0.5mm, color = red, preaction={-triangle 60,thin,draw,shorten >=-0.5mm,color=red}] ($(a)+(1.6,0)$) -- ($(a)+(1.6,-0.5)$);
	\draw[line width = 0.5mm, color = red, preaction={-triangle 60,thin,draw,shorten >=-0.5mm,color=red}] ($(a)+(2,0)$) -- ($(a)+(2,-0.5)$);
	\draw[line width = 0.5mm, color = red, preaction={-triangle 60,thin,draw,shorten >=-0.5mm,color=red}] ($(a)+(2.4,0)$) -- ($(a)+(2.4,-0.5)$);
	\draw[line width = 0.5mm, color = red, preaction={-triangle 60,thin,draw,shorten >=-0.5mm,color=red}] ($(a)+(2.8,0)$) -- ($(a)+(2.8,-0.5)$);
	\draw[line width = 0.5mm, color = red, preaction={-triangle 60,thin,draw,shorten >=-0.5mm,color=red}] ($(a)+(3.2,0)$) -- ($(a)+(3.2,-0.5)$);
	\draw[line width = 0.5mm, color = red, preaction={-triangle 60,thin,draw,shorten >=-0.5mm,color=red}] ($(a)+(3.6,0)$) -- ($(a)+(3.6,-0.5)$);
	\draw[line width = 0.5mm, color = red, preaction={-triangle 60,thin,draw,shorten >=-0.5mm,color=red}] ($(a)+(4,0)$) -- ($(a)+(4,-0.5)$);
	\draw[line width = 0.5mm, preaction={-triangle 60,thin,draw,shorten >=-0.5mm}] ($(b)+(-4.95,-1.05)$) -- ($(b)+(-4.3,-0.4) $);
	\coordinate (c) at ($(b)+(-4.95,-1.35)$);
	\coordinate (d) at ($(b)+(4.95,-1.35)$);
	\node at (c) {adiabatic};
	\node at (b|-c) {temperature (output $y$)};
	\node at (d) {$T=0K$};
	\draw[line width = 0.5mm, preaction={-triangle 60,thin,draw,shorten >=-0.5mm}] ($(b)+(0,-1.05)$) -- ($(b)+(0,-0.4) $);
	\draw[line width = 0.5mm, preaction={-triangle 60,thin,draw,shorten >=-0.5mm}] ($(b)+(4.95,-1.05)$) -- ($(b)+(4.3,-0.4) $);
	\end{tikzpicture}
\caption{1D-Beam with heat source (input) and measured temperature (output)}
\label{fig:HeatTransfer}
\end{figure}

In this section the Loewner framework, endowed with the index reduction and regularization procedure outlined in Section \ref{sec:Matlab}, is illustrated by means of an example from the Oberwolfach Model Reduction Benchmark Collection \cite{BenMS05}. We consider the nonlinear heat transfer in a one-dimensional beam discussed in \cite{LieYK05}. A schematic illustration of the system is depicted in Figure \ref{fig:HeatTransfer}. For the sake of simplicity we restrict ourselves to the single-input single-output (SISO) case in contrast to the multiple-input multiple-output (MIMO) system considered in \cite{LieYK05}.

The governing equation of the physical system is a parabolic partial differential equation describing the temporal progress of the spatial temperature distribution along the beam. However, instead of the absolute temperature $T_{abs}$, a relative temperature $T$ is considered, i.\,e., $T = T_{abs}-T_{ref}$ with respect to a reference temperature $T_{ref}$. The initial condition is chosen homogeneously as $T=0$ over the whole beam at time $t_0=0$. Furthermore, at the left boundary an adiabatic end is assumed, i.\,e., zero temperature gradient, and at the right boundary the relative temperature is equal to zero for all times $t>0$ \cite{LieYK05}.

In this example we are rather interested in the input-output (I/O) behavior than in the time progress of the entire temperature distribution. As an input a heat source is applied affecting the whole beam homogeneously and the temperature at the middle of the beam represents the system output.

Moreover, a nonlinearity comes into play by considering a thermal conductivity which depends on the temperature polynomially, i.\,e.,

\begin{equation*}
k\left(T\right) = \sum\limits^{N}_{i=0}a_iT^i
\end{equation*}
with given coefficients $a_i$. After modeling, discretization and renaming of variables ($T\rightarrow x$) one obtains a dynamical system of the form

\begin{equation}
\begin{aligned}
	E\dot{x} &= Ax+bu+f\left(x\right),\\
	y &= c^Tx,
\end{aligned}
\label{eq:controlSystem}
\end{equation}
where $x\in\mathbb{R}^{n}$ denotes the state vector (discrete approximation of temperature), $u$ the input (heat source) and $y$ the output (temperature at the middle of the beam). Furthermore, $E,A\in\mathbb{R}^{n,n}$ and $b,c\in\mathbb{R}^n$ represent the linear part and the function $f:\mathbb{R}^n\rightarrow\mathbb{R}^n$ constitutes the nonlinear part of the dynamical system. More details regarding the modeling and discretization may be found in \cite{LieYK05}.

Depending on the mesh size, there are two systems of different dimensions available within the Oberwolfach Model Reduction Benchmark Collection: $n=15$ and $n=410$. Since the main intention is to illustrate the need for the regularized Loewner approach, we choose the system of dimension $n=15$ due to the significantly smaller simulation times.

In order to use the Loewner method we need sampled data of the transfer function of the system. Since a nonlinear dynamical system is considered, there is only little hope to find an analytic expression for the transfer function of the system. The idea is instead to utilize the system's impulse response and determine a linear transfer function describing the input-output behavior of the system for the chosen input. Due to the nonlinearity of the system the obtained transfer function has only a limited validity range with its size depending on the impact of the nonlinearity on the I/O map.

Since an actual impulse response is numerically unfeasible, instead we create the step response and differentiate it numerically, in order to obtain an approximation of the impulse response, as in \cite{BorCG12}. The discrete values of the impulse response are equal to the Markov parameters $h_k$ ($k=0,1,\ldots$) of the corresponding discrete-time system leading to the discrete-time transfer function

\begin{equation}
\tilde{H}\left(z\right) = \sum\limits_{k=0}^{\infty}h_kz^{-k}.
\label{eq:Markov}
\end{equation}
Since the impulse response of the considered system approaches zero for large time values, the same holds for the Markov parameters with high index. Consequently, the infinite sum of equation (\ref{eq:Markov}) may be truncated while retaining a reasonable level of accuracy. For applying the Loewner approach, the transfer function is expected to map from the Laplace transforms of the inputs to the Laplace transforms of the outputs, cf. Section \ref{intro}. However, the obtained transfer function $\tilde{H}$ refers to the Z-domain. In order to obtain an expression for the transfer function of the continuous-time system the bilinear transformation is used to transform from the Z-domain to the Laplace domain \cite{OppSB99}, i.\,e.,
\begin{equation*}
z = \frac{1+\frac{\Delta t}{2}s}{1-\frac{\Delta t}{2}s}
\end{equation*}
leading to
\begin{equation*}
H\left(s\right) = \tilde{H}\left(\frac{1+\frac{\Delta t}{2}s}{1-\frac{\Delta t}{2}s}\right),
\end{equation*}
where $\Delta t$ denotes the sampling time interval, which is equal to the time step size used for the simulation of the step response.
After these preliminary steps one obtains an approximate transfer function $H\left(s\right)$ which may be sampled in order to apply the Loewner framework as described in Section \ref{sec:realization}.

The aforementioned procedure to determine a linear approximation of the transfer function is based on simulating the step response, i.\,e., using the Heaviside step function $\Theta\left(t\right)$ as input. However, we prefer to consider multiples of the Heaviside function as in \cite{LieYK05}. For this purpose, we split the input $u\left(t\right) = a\Theta\left(t\right)$ by putting the constant factor $a$ into the $b$-vector, leading to a system equivalent to (\ref{eq:controlSystem}), but replacing $b$ by $\tilde{b}=ab$ and $u\left(t\right)$ by $\tilde{u}\left(t\right)=\Theta\left(t\right)$. Consequently, we may consider the  step response without being restricted to an input magnitude of $ 1\, W$.

\begin{figure}
\centering
\includegraphics[scale=0.28]{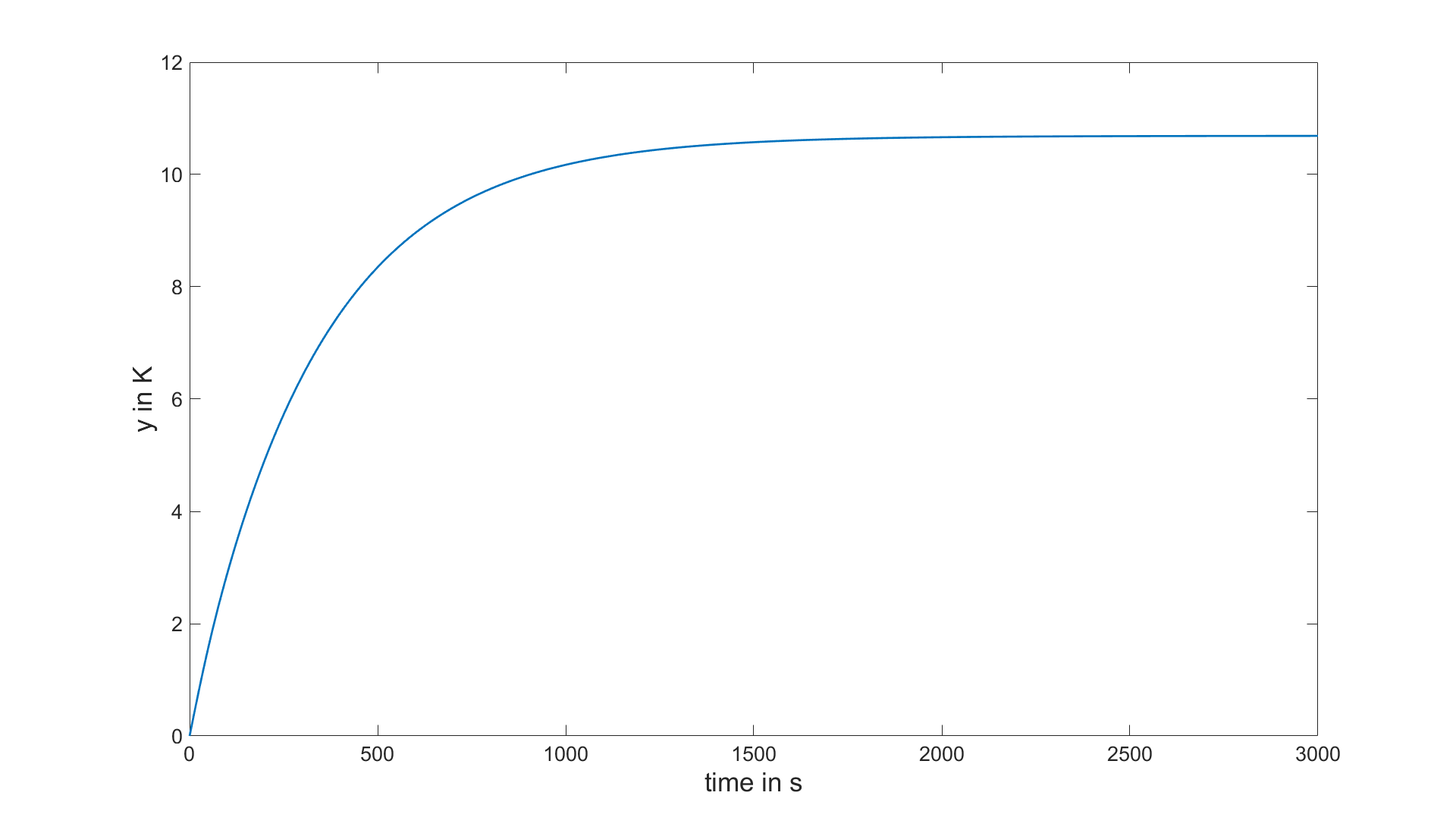}
\caption{Step response of the full model ($u=10^5\, W$)}
\label{fig:stepResponse_fullModel}
\end{figure}

\begin{figure}
\centering
\includegraphics[scale=0.27]{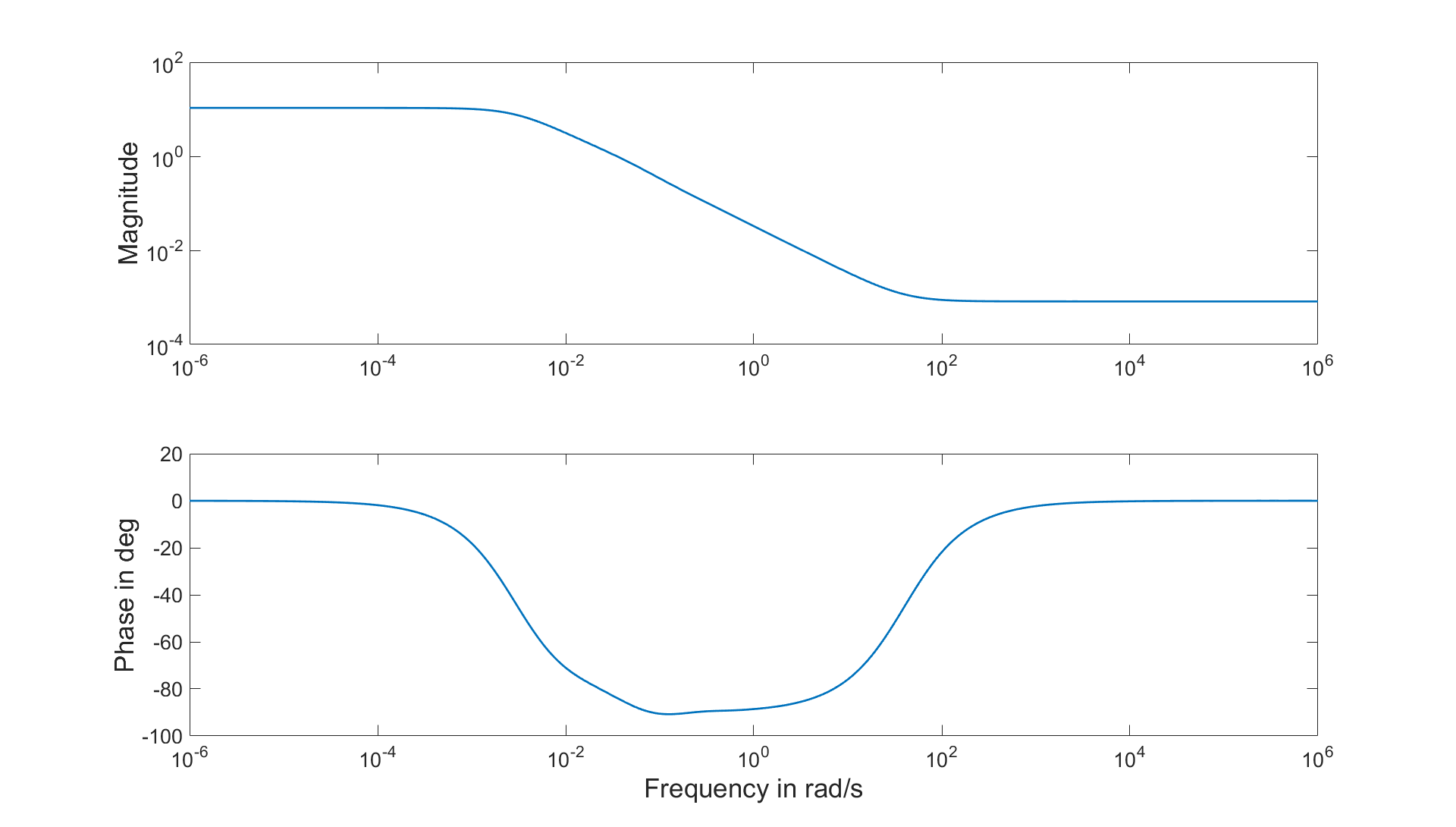}
\caption{Bode plot of approximate transfer function}
\label{fig:BodePlot_transferFunction}
\end{figure}

We determine the step response using an input step of $10^5\, W$ which is within the range of heat source magnitudes considered in \cite{LieYK05}. The corresponding output step response is given in Figure \ref{fig:stepResponse_fullModel}. Applying the procedure outlined above, we obtain an approximation of the transfer function based on this step response. The Bode plot of this transfer function is depicted in Figure \ref{fig:BodePlot_transferFunction}.

Based on the linear approximation of the transfer function, we apply the Loewner framework to obtain a low-dimensional realization which interpolates the transfer function. For this, we only need to choose interpolation points but no tangential interpolation directions, since we only have one input and one output (SISO case). As we would like to approximate the transfer function over a wide range of frequencies, logarithmically equidistant sets of interpolation points are chosen. Moreover, in order to be able to check the interpolation easily by means of the Bode plot, purely imaginary numbers are chosen for the interpolation points. Furthermore, the complex-conjugate interpolation points are added automatically and a coordinate transformation is performed (in function \verb+loewner_mod+) to obtain a real-valued realization. Numerous constellations of interpolation point sets have been tested in an automatic fashion to get a better insight into proper selections of interpolation point ranges.

Comprehensive tests have shown that the range of interpolation points should not be chosen wider than $20$ orders of magnitude. Ranges  that are too wide lead numerically to a violation of the rank conditions which are necessary for the Loewner approach to be applicable, cf. \eqref{eq:rankCondition}. In accordance with this observation, the following rule of thumb may be formulated: The smaller the range of interpolation points, the higher the admissible number of interpolation points.

\begin{figure}
\centering
\includegraphics[scale=0.28]{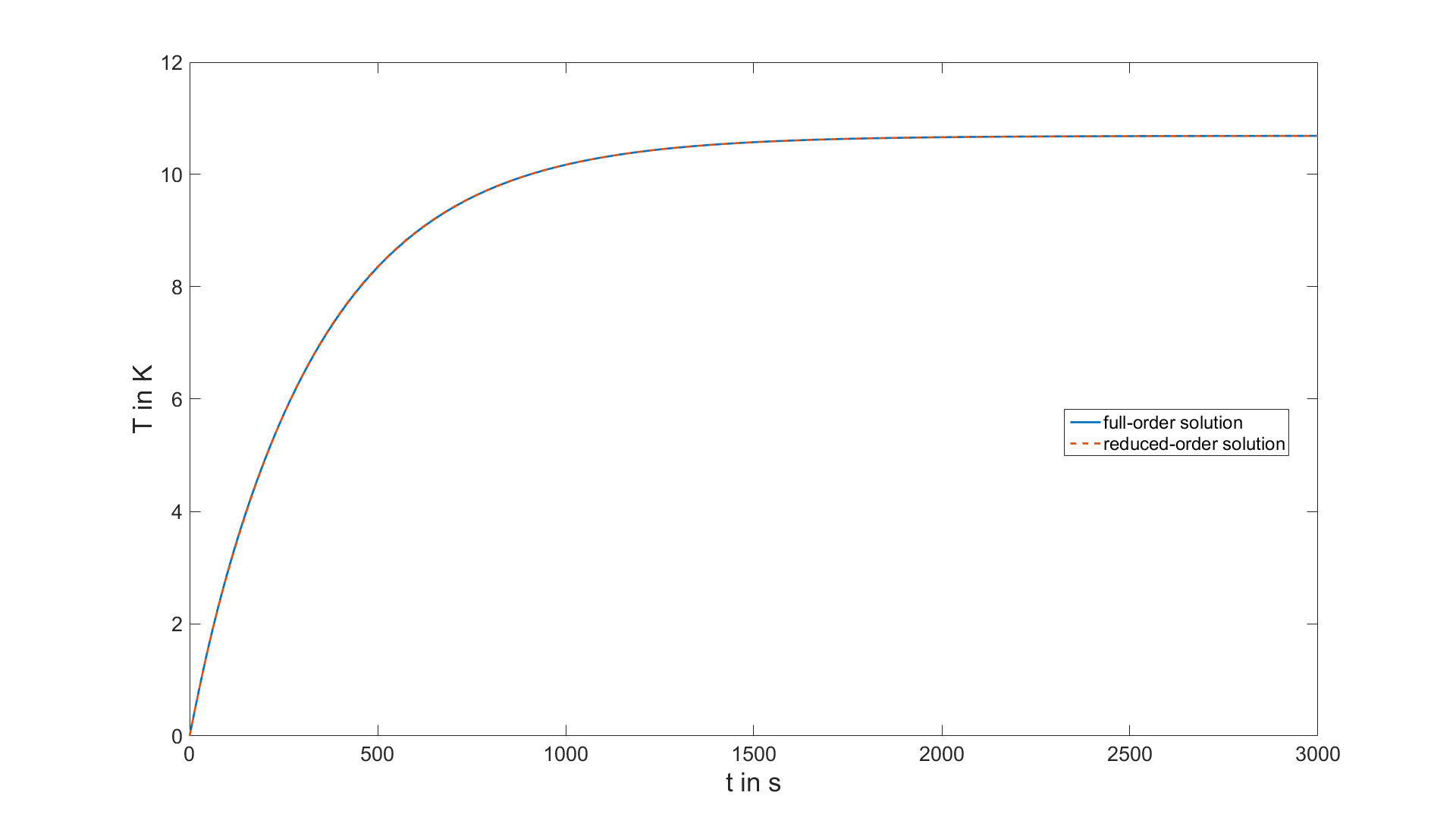}
\caption{Comparison of step responses ($u=10^5\, W$)}
\label{fig:stepResponse_comparison}
\end{figure}

In addition to this, it should be noted that the number of interpolation points is proportional to the dimension of the Loewner realization, at least in the regular case (cf. Section \ref{sec:realizationFunc}). Therefore, we are mainly interested in interpolation point sets containing only a small number of points. In order to obtain real-valued realization matrices, four is the minimal number of interpolation points needed. Several constellations have been tested. The smallest step response error (measured in the maximum norm) is provided by the interpolation point set
\begin{equation}
S=\left\lbrace -10^{3.6}i;-10^{1.8}i;-10^{-1.8}i;-10^{-3.6}i;10^{-3.6}i;10^{-1.8}i;10^{1.8}i;10^{3.6}i\right\rbrace.
\label{eq:optimalPointSet}
\end{equation}
The dimension of the Loewner realization is half the number of interpolation points which leads to a state space dimension of four in this case.
The comparison of the step response of the reduced system to that of the original system is presented in Figure \ref{fig:stepResponse_comparison}. It should be noted that the step height for the original model is again $u=10^5\, W$, whereas the step height for the reduced model is $u=1$, cf. discussion above about multiples of the Heaviside function. The excellent agreement of the step responses is obvious. We emphasize that a nonlinear system of dimension $15$ has been reduced to a linear system of dimension four.

It is noteworthy that the reduced system only provides a good approximation for the I/O behavior of the full system, whereas the internal state variables of the original model are not captured in the reduced order model. However, in many applications, approximating the I/O behavior is sufficient, e.\,g., in control applications.

For the case of the interpolation point set \eqref{eq:optimalPointSet}, the Loewner realization is strangeness-free as well as completely controllable and observable. Thus, the index reduction and regularization procedure is not necessary in this case. In contrast, some interpolation point sets lead to realizations with strangeness-index greater than or equal to one. One example set is given by the set
\begin{equation*}
\Sigma=\left\lbrace -10^{7}i;-10^{-7}i;10^{-7}i;10^{7}i\right\rbrace.
\end{equation*}
Without the index reduction procedure, the numerical integration of the resulting Loewner realization by means of the {\tt Matlab} solver \verb+ode15s+ fails due to the higher index. However, the regularization procedure transforms the reduced system to an equivalent strangeness-free system and the numerical integration succeeds. This example emphasizes the need of a regularization procedure.

\begin{figure}
\centering
\includegraphics[scale=0.28]{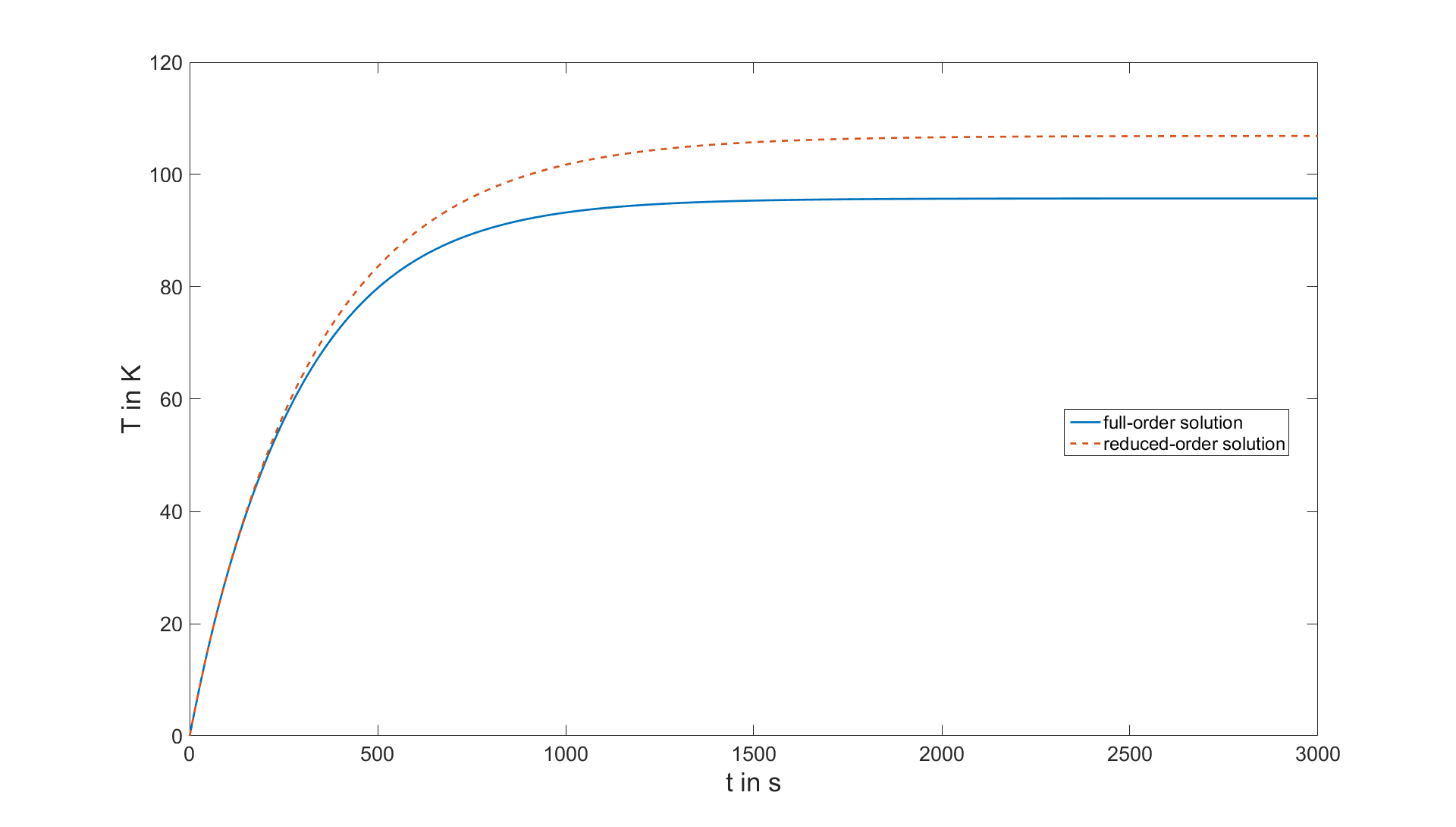}
\caption{Comparison of step responses ($u=10^6\, W$)}
\label{fig:stepResponse_comparison_u10}
\end{figure}

On top of potential higher-index, often unstable realizations are obtained, which are not avoided by the regularization procedure presented in this work. These systems lead to trouble when simulating the step response due to the unstable behavior. This directly leads to the topic of stability-preserving model reduction. This is not within the scope of this report but for completeness we mention the passivity-preserving interpolation approach in \cite{Ant05}. It also preserves stability and is based on choosing the spectral zeros of the original transfer function as interpolation points. The spectral zeros are defined as the solutions of the equation
\begin{equation*}
H\left(s\right)+H\left(-s\right) = 0.
\end{equation*}
As a final remark of this section, it should be emphasized that the determined transfer function and the resulting reduced order model are only valid for inputs being close to the test input $u=10^5\, W$. When considering much bigger or smaller input steps, the difference between the step responses of the reduced and the full system are significantly larger. The reason for this is the nonlinearity of the original model, which can be approximated by a linear model only locally. To illustrate this discrepancy, Figure \ref{fig:stepResponse_comparison_u10} shows the comparison of the step responses for an input step of $u=10^6\, W$ where the reduced model is the same as in Figure \ref{fig:stepResponse_comparison} (based on step response with $u=10^5\, W$). The qualitative behavior is indeed well approximated by the reduced model but the quantitative agreement is bad when the height of the input step is much larger ($u=10^6\, W$) than that used for determining the reduced model ($u=10^5\, W$).
In order to approximate the original system for a wide range of inputs, several linear surrogate models are needed or an approach different from the basic Loewner framework has to be applied.

\section{Conclusion}

In order to make the realization obtained from the Loewner framework suitable for simulation and control applications, we have presented a regularization procedure resulting in a strangeness-free as well as completely controllable and observable system. This procedure has been implemented in \texttt{Matlab} and is illustrated by means of a nonlinear heat transfer problem. The numerical results reveal that applying the pure Loewner realization may lead to higher-index or not completely controllable or observable systems. When using the regularization procedure presented in this work, the realization is transformed to an equivalent system being strangeness-free and completely controllable and observable. These properties are important when performing simulations or when applying control methods based on the Loewner realization.\\

\noindent\textbf{Acknowledgements.} The authors gratefully acknowledge the support by the Deutsche Forschungsgemeinschaft (DFG) as part of the collaborative research center SFB 1029 \textit{Substantial efficiency increase in gas turbines through direct use of coupled unsteady combustion and flow dynamics}, project A02 \textit{Development of a reduced order model of pulsed detonation combuster}.

\bibliographystyle{plain}

\end{document}